\documentclass{ifacconf}

\usepackage{graphicx}      
\graphicspath{{./figs/}}
\usepackage{natbib}        
\usepackage{xcolor}
\usepackage[utf8]{inputenc}

\usepackage{fancyhdr}
\usepackage{subfigure}
\makeatletter
\let\old@ssect\@ssect 
\makeatother
\usepackage{hyperref}
\makeatletter
\def\@ssect#1#2#3#4#5#6{%
  \NR@gettitle{#6}
  \old@ssect{#1}{#2}{#3}{#4}{#5}{#6}
}
\makeatother

\usepackage{amsmath,amssymb,amsfonts}
\usepackage{algorithmic}
\usepackage{textcomp}
\newcommand{\vect}[1]{\boldsymbol{\mathbf{#1}}}
\usepackage{mathtools}
\usepackage{placeins}
\usepackage{siunitx}

\newcommand{\ui}[2]{#1_{\mathrm{#2}}}

\begin{document}
\begin{frontmatter}

    \title{Model Predictive Control of a Vehicle using Koopman Operator}

    \thanks[footnoteinfo]{© 2020 the authors. This work has been accepted to IFAC for publication under a Creative Commons Licence CC-BY-NC-ND}        

    \author[1]{Vít Cibulka}
    \author[1]{Tomáš Haniš}
    \author[1,3]{Milan Korda}
    \author[1]{Martin Hromčík}
    \address[1]{Dept. of Control Engineering, Faculty of Electrical Engineering, \\Czech Technical University in Prague, The Czech Republic \\(emails: cibulka.vit@fel.cvut.cz, 
    hanis.tomas@fel.cvut.cz, korda.milan@fel.cvut.cz, hromcik.martin@fel.cvut.cz)}
    \address[3]{CNRS, Laboratory for Analysis and Architecture of Systems, Toulouse, France \\(email: korda@laas.fr)}

    \begin{abstract}                
        This paper continues in the work from \cite{cibulka_koop_identification} where a nonlinear vehicle model was approximated
        in a purely data-driven manner by a linear predictor of higher order, namely the Koopman operator.
        The vehicle system typically features a lot of nonlinearities such as rigid-body dynamics,
        coordinate system transformations and most importantly the tire.

        These nonlinearities are approximated in a predefined subset of the state-space by the \textit{linear} Koopman operator and
        used for a \textit{linear} Model Predictive Control (MPC) design in the high-dimension state space where the nonlinear system dynamics
        evolve \textit{linearly}. The result is a nonlinear MPC designed by linear methodologies.

        It is demonstrated that the Koopman-based controller is able to recover from a very unusual state of the vehicle where all the
        aforementioned nonlinearities are dominant.
        The controller is compared with a controller based on a classic local linearization and shortcomings of this
        approach are discussed.
    \end{abstract}


\end{frontmatter}

\section{Introduction}
A vehicle is a nonlinear system that is becoming more interesting from the control engineering point of view with the
ever increasing number of electric vehicles.
This gives an opportunity for sophisticated control systems to take the place of old-fashioned solutions
which are currently present in the majority of vehicles today.

This paper examines the nonlinear control of the vehicle described by a linear predictor which is valid in a predefined subset state space,
which allows for exploitation of linear control methods on the nonlinear system.
The linear predictor used in this paper is the Koopman operator (\cite{Koopman1931}).

The Koopman operator, an
increasingly popular tool for global linearization and analysis of nonlinear dynamics (\cite{Mezic2005}, \cite{Korda_opt_g} ,\cite{Korda_koop}, \cite{Mezic_control_suggestion}), is used in this work
to approximate the vehicle nonlinear dynamics in order to achieve
a linear representation of the system in a predefined subspace of the state space.

This paper continues in the work from \cite{cibulka_koop_identification}, where different methods for global
linearization of the single-track model were used. The most promising method (described in detail in \cite{Korda_opt_g})
is used for approximation of autonomous and controlled behaviour of the nonlinear vehicle system by a high-dimensional linear
system.
The resulting linear system is then used for \textit{linear} Model Predictive Control (MPC) design and verified
against a MPC based on local linearization which was the prevalent approach of tackling nonlinear
systems in the past.

\section{Single-track model}
\label{sec:singletrack_definition}
\begin{figure}[htbp]
    \centerline{\includegraphics[width=0.5\textwidth]{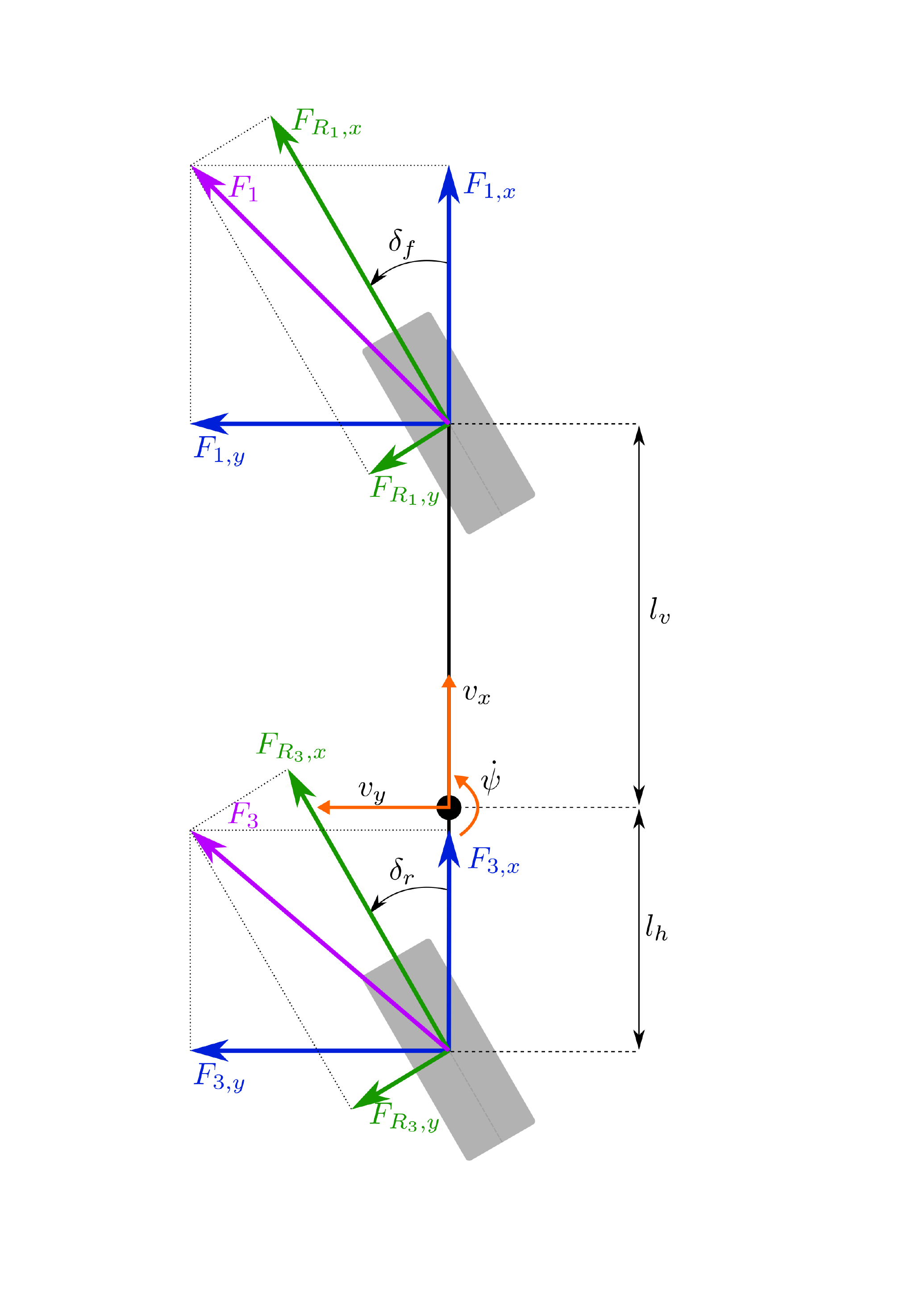}}
    \caption{The single-track model. Forces $F_{R_2}$ and $F_{R_4}$ are not depicted in the figure because in a general case with symmetric tires $F_{R_2} = F_{R_1}$ and $F_{R_4} = F_{R_3}$.}
    \label{fig:singletrack}
\end{figure}

The vehicle model derived in \cite{cibulka_koop_identification} will be reviewed here.
The model is depicted in Fig.~\ref{fig:singletrack}.
State vector of the single-track model is 
\begin{equation}
    \begin{bmatrix}
        v_x (\SI{}{\metre\per\second})          , &
        v_y (\SI{}{\meter\per\second})          , &
        \dot{\psi}_z (\SI{}{\radian\per\second})
    \end{bmatrix}^\top,
\end{equation}
 where $\ui{v}{x}$ is longitudinal velocity, $\ui{v}{y}$ lateral velocity and $\ui{\dot{\psi}}{z}$ is 
  yawrate.
Inputs to the model are rear longitudinal slip ratios \(\kappa_\mathrm{ r }\) and front steering angle $\delta_{\mathrm{ f }}$.

The vehicle body is modeled as a rigid body using Newton-Euler equations
\begin{equation}
    m_v (
    \begin{bmatrix}
        \ui{ \dot{v} }{ x } \\
        \dot{v}_\mathrm{ y } \\
    \end{bmatrix} +
    \dot{\psi}_\mathrm{ z }
    \begin{bmatrix}
        -v_\mathrm{ y } \\
        v_\mathrm{ x }  \\
    \end{bmatrix} )
    =
    \sum_{i=1}^{4}
    \begin{bmatrix}
        {F}_{i,\mathrm{ x }} \\
        {F}_{i,\mathrm{ y }} \\
    \end{bmatrix}
    -\frac{1}{2}c_\mathrm{ w } \rho A_\mathrm{ w } \sqrt{v_\mathrm{ x }^2 + v_\mathrm{ y }^2}
    \begin{bmatrix}
        v_\mathrm{ x } \\v_\mathrm{ y }
    \end{bmatrix} 
\end{equation} and
\begin{equation}
    J_\mathrm{zz} \ddot{\psi}_\mathrm{ z }
    =\sum_{i=1}^{4} \vect{r}_i \mathbf{F}_{i}
    \label{eq:yawr},
\end{equation}
where
\begin{equation}
    \vect{r} =
    \begin{bmatrix}
        \vect{r_1} & \vect{r_2} & \vect{r_3}& \vect{r_4}
    \end{bmatrix}
    =
    \begin{bmatrix}
        \begin{bmatrix}
            l_\mathrm{ v } \\0\\0
        \end{bmatrix},
        \begin{bmatrix}
            l_\mathrm{ v } \\0\\0
        \end{bmatrix},
        \begin{bmatrix}
            -l_\mathrm{ h } \\0\\0
        \end{bmatrix},
        \begin{bmatrix}
            -l_\mathrm{ h } \\0\\0
        \end{bmatrix}
    \end{bmatrix}
\end{equation}

is the vector describing position of each wheel with respect to the center of gravity and
\(\vect{F}_i = \begin{bmatrix}
    F_{i,x}\\F_{i,y}
\end{bmatrix}\) is a vector of forces acting on \(i^{\mathrm{th}}\) wheel. The vector and its elements are depicted in Fig.~\ref{fig:singletrack}.
Note that although Fig.~\ref{fig:singletrack} might suggest 
that the model has 2 wheels, it is defined with 4 wheels, where the left and right wheels are in the same place.
This allows for usage of asymmetrical tire models (such as the one used in this paper).
The parameters $l_\mathrm{ v }$ and $l_\mathrm{ h }$ are distances of wheels from CG, as depicted in Fig.~\ref{fig:singletrack}. The wheels are numbered in this order: front-left, front-right, rear-left, rear-right.
$ m_\mathrm{ v } $ is the vehicle mass,
\(\mathbf{F}_{i,\mathrm{x/y}}\) is a force acting on i-th wheel along x/y axis in body-fixed coordinates. $F_{R_i,\mathrm{ x }}$ is a force acting along x axis in wheel coordinate system.
The term \( -\frac{1}{2}c_\mathrm{ w } \rho A_\mathrm{ w } \sqrt{v_\mathrm{ x }^2 + v_\mathrm{ y }^2}
\begin{bmatrix}
    v_\mathrm{ x } \\v_\mathrm{ y }
\end{bmatrix}\)
is an approximation of air-resistance, \(c_\mathrm{ w }\) is a drag coefficient, \(\rho\) is air density and \(A_\mathrm{ w }\) is the total surface exposed to the air flow.
\(J_\mathrm{zz}\) is the vehicle inertia about z-axis and
\(J_{R_i}\) is the wheel inertia about y-axis.

The forces \(  \begin{bmatrix}
    F_{R_i,\mathrm{ x }} \\
    F_{R_i,\mathrm{ y }} \\
\end{bmatrix} \)
are calculated using the ``Pacejka magic formula'' \cite{Pacejka2012}
\begin{equation}
    F = D \cos(C \arctan(Bx - E(Bx - \arctan(Bx)))).
    \label{pacejka}\end{equation}
The same formula can be used for calculating $F_{R_i,\mathrm{ x }}$ (tire longitudinal force) and $F_{R_i,\mathrm{ y }}$ (tire lateral force) with a different set of parameters for each.
The argument $x$ can be either sideslip angle $\alpha$ or longitudinal slip ratio $\kappa$ (usually denoted as \(\lambda\) which is used for eigenvalue in this paper) (see \cite{Pacejka2012}) for calculating $F_{R_i,\mathrm{ y }}$ or $F_\mathrm{ x }$ respectively.
The parameters $B,C,D$ and $E$ are generally time-dependent. This work uses the Pacejka tire model \cite{Pacejka2012} with coefficients from the \textit{Automotive challenge 2018} organized by Rimac Automobili.
The transformation of tire forces from wheel-coordinate system to car coordinate system is done as follows
\begin{equation}
    \begin{bmatrix}
        F_{i,\mathrm{ x }} \\
        F_{i,\mathrm{ y }} \\
    \end{bmatrix}
    =
    \begin{bmatrix}
        \cos(\delta_i) & -\sin(\delta_i) \\
        \sin(\delta_i) & \cos(\delta_i)
    \end{bmatrix}
    \begin{bmatrix}
        F_{R_i,\mathrm{ x }} \\
        F_{R_i,\mathrm{ y }} \\
    \end{bmatrix}.
\end{equation}

\section{Linear predictors}
Linear predictor is a linear model of a controlled system that is able to provide the prediction of the future 
behaviour of the controlled system with sufficient accuracy.
The predictor used in this paper is the Koopman operator and it will be used as a control design model 
for MPC. The Koopman operator is infinite-dimensional linear system, which is able to describe the nonlinear behaviour 
of the controlled system. A finite-dimensional approximation of the Koopman operator will be used as a control design model 
for a linear MPC resulting in a control law that is linear in the state space of the Koopman operator, but nonlinear 
in the original state space of the nonlinear controlled system.

\subsection{Koopman operator}
The Koopman operator is used as a linear predictor of the nonlinear dynamics of the system Sec.~\ref{sec:singletrack_definition}.
The basic idea consists in transforming (\textit{lifting}) the nonlinear state space to a new high-dimensional, linearly evolving state space.
The control design is then performed in the linear state space using linear control methodology.
\label{sec:eigenfun}
Let us assume a discrete nonlinear uncontrolled system with state
\(x_k\) at time step \(k\), dynamics \(f_{\mathrm{u}}(.)\), output \(y_k\) and
output equation \(h( {x_k} )\):
\begin{align}
    \begin{split}
        \label{eq:f_u}
        x_{k+1} &= f_{\mathrm{u}}(x_k) \\
        y_k &= h(x_k).
    \end{split}
\end{align}

The Koopman operator \(\mathcal{K}: \mathcal{C}(\mathbb{R}^n) \rightarrow  \mathcal{C}(\mathbb{R}^n)\),
with \(\mathcal{C}({\mathbb{R}^n})\) denoting a space of continuous functions defined on \(\mathbb{R}^n\), is defined as
\begin{equation}
    \label{eq:koop_def}
    ( \mathcal{K}\phi )(x_k) = \phi(f_{\rm{u}}(x_k))
\end{equation}
for each basis function \(\phi: \mathbb{R}^n \rightarrow \mathbb{R}\) where \(n\) is size of the state vector \(x_k\).
In our case, the function \(\phi\) will also be an \textit{eigenfunction} of the
operator \(\mathcal{K}\), meaning that the following holds:
\begin{equation}
    \label{eq:phi_linear_evolution}
    \phi(x_{k+1}) = \lambda \phi(x_k),
\end{equation}
for some eigenvalue \(\lambda \in \mathcal{R}\).
The functions \(\phi\) will be constructed from trajectories of \eqref{eq:f_u}
according to
\begin{equation}
    \label{eq:basis}
    \phi(x^j_k) = \phi(x^j_k)_{\lambda,g} = \lambda^{k}g_\phi(x^j_0),
\end{equation}
where \(j\) is a trajectory of \eqref{eq:f_u} starting in $x^j_0$ and \(x_k^j\) is the point to which the system will get after \(k\) time-steps.
The state vector $x^j_k$ is transformed with a function $\phi(x^j_k)$, defined
according to \eqref{eq:basis}
for an arbitrary eigenvalue $\lambda$ and an arbitrary function $g_\phi : \mathbb{R}^n \rightarrow \mathbb{R}$.

Note that the definition from \eqref{eq:basis} fulfills the requirement of \eqref{eq:phi_linear_evolution} because
\begin{equation}
    \phi(x_{k+1}^j) =  \lambda^{k+1} g_\phi(x_0^j) = \lambda \cdot \lambda^{k} g_\phi(x_0^j)  = \lambda \phi(x_k^j).
\end{equation}
In other words, $\phi(x_k^j)$ evolves linearly along trajectories of the system \eqref{eq:f_u}.
The trajectories must fulfill  certain assumptions in order for the definition \eqref{eq:basis} 
to be valid. The assumptions are beyond the scope of this paper and are discussed in \cite{Korda_opt_g}.

\subsection{Uncontrolled case}
\label{sec:uncontrolled}
The functions \(g_{\phi}(.)\) can be replaced with scalars because they are evaluated only at the starting points \(x_0^{j}\) of the trajectories \(j\).
Let us denote the set of starting points \(x_0^j\) as \(\Gamma\). The evaluation of \(g_{\phi}(.)\) on a point from \(\Gamma\)
will be denoted as
\begin{equation}
    g_{p,i}^{j} = g_{\phi}(x_0^j), \text{ for } x_0^j \in \Gamma,
\end{equation}
where \(p\) denotes the number of output \((p = 1,2,...N_y)\) with \(N_y\) being the total number of outputs and \(i\) is the associated eigenvalue.
The association of \(g_{p,i}^{j}\) with a specific eigenvalue and a specific output allows for a trivial derivation of the \(A\) and \(C\) matrices,
which will be discussed further below.
The values \(g_{p,i}^{j}\) can be optimized in a \textit{convex} manner in order to approximate the output values by
\begin{equation}
    \label{eq:g_output_eq}
    y_{p,k}^{j} = \sum_{i=1}^{N_\Lambda} \lambda_i^k g_{p,i}^j,
\end{equation}
where \( y_{p,k}^{j}\) is the \(p^{\mathrm{th}}\) output of \(j^{th}\) trajectory at time-step \(k\) and \(N_\Lambda\) is the number of eigenvalues.
The solution of \eqref{eq:g_output_eq} for output \(p\) can be written in matrix form as
\begin{equation}
    \label{eq:ch3:opt_g}
    || L g_{p} - F_p ||_2^2 + \zeta||g_p||_2^2,
\end{equation}
where \(L\) is a matrix containing the eigenvalues \(\lambda\), \(F_p\) is a matrix
of outputs from all trajectories and \(\zeta\) is a regularization term.
The optimized value is the vector \(g_p\) which contains \(g_{p,i}^{j}\) for all \(\lambda_i\) and all trajectories.

The concrete form of the matrices in \eqref{eq:ch3:opt_g} 
 can be found in \cite{Korda_opt_g}, as well as the algorithm 
 for finding the eigenvalues \(\lambda_i\).

In order to obtain the matrices A and C consider the eigenfunction definition from \eqref{eq:phi_linear_evolution}, for 
\begin{equation}
\boldsymbol{\phi}(x^j_k)~=~\begin{bmatrix}
        \phi_1(x^j_k) & \phi_2(x^j_k) & \ldots \phi_{N_\phi}(x^j_k)
    \end{bmatrix}^\top
\end{equation}
the dynamics
\begin{equation}
    z_{k+1} = A z_k
\end{equation}
can then be written as
\begin{equation}
    \begin{bmatrix}
        \phi_1(x^j_{k+1}) \\
        \phi_2(x^j_{k+1}) \\
        \vdots            \\
        \phi_{N_\phi}(x^j_{k+1})
    \end{bmatrix}
    =
    \begin{bmatrix}
        \lambda_1                                \\
         & \lambda_2                             \\
         &           & \ddots                    \\
         &           &        & \lambda_{N_\phi}
    \end{bmatrix}
    \begin{bmatrix}
        \phi_1(x^j_k) \\
        \phi_2(x^j_k) \\
        \vdots        \\
        \phi_{N_\phi}(x^j_k)
    \end{bmatrix}.
\end{equation}

Choosing \eqref{eq:phi_linear_evolution} as basis functions immediately yields the diagonal $A$ matrix.
The output matrix C is also trivial thanks to \eqref{eq:g_output_eq}.
\begin{equation}
    C =
    \begin{bmatrix}
        1\dots 1 &          &          \\
                 & 1\dots 1 &          \\
                 &          & 1\dots 1
    \end{bmatrix}_{N_y \times (N_y \cdot N_\Lambda)}.
\end{equation}
Note that in this case, the Koopman operator defined in \eqref{eq:koop_def} is implemented as the state matrix \(A\).

\subsection{Controlled case}
\label{sec:koop_controlled}
\begin{figure}[htbp]
    \centerline{\includegraphics[width=0.5\textwidth]{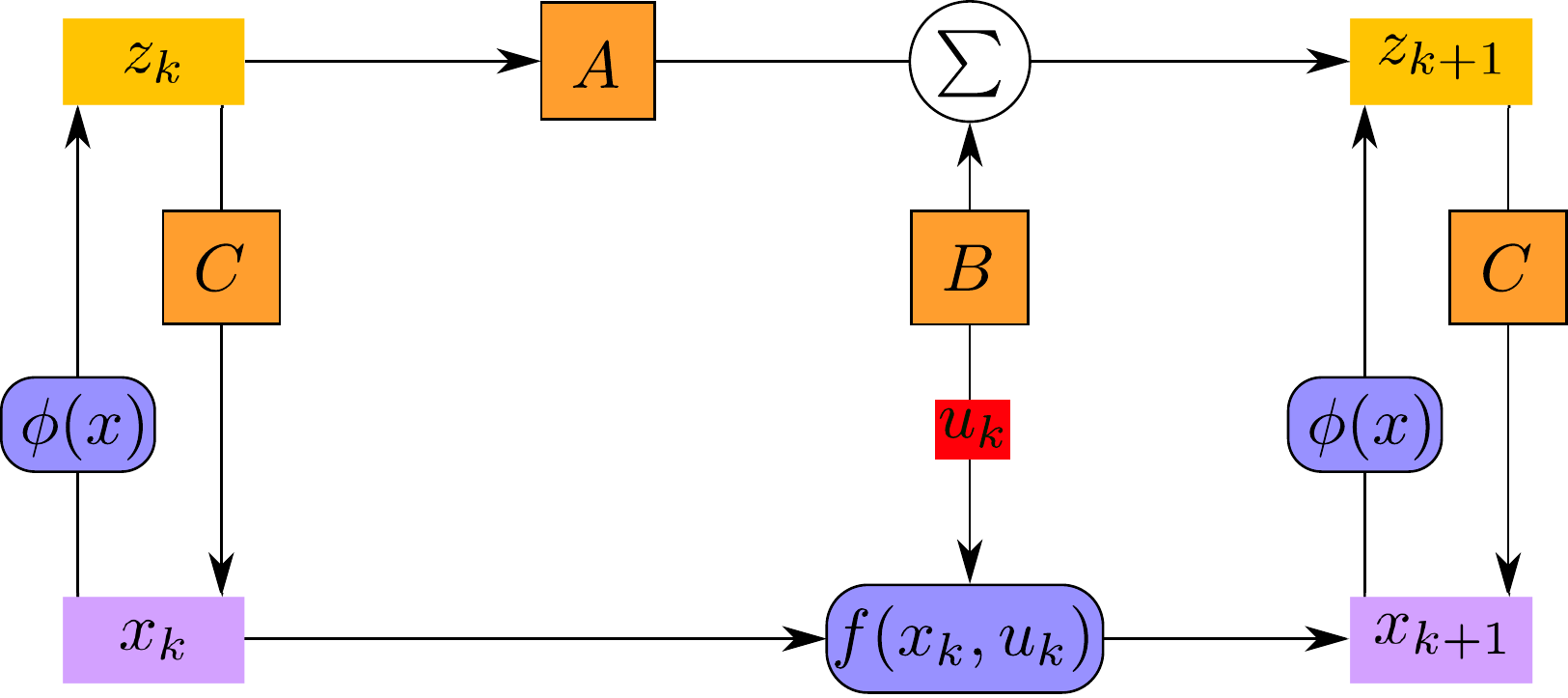}}
    \caption{Discrete-time scheme showing the relationship of a nonlinear system and
        its linear approximation. }
    \label{fig:koopman}
\end{figure}

In this work however, a controlled scenario will be considered.
The discrete \textit{ nonlinear } controlled system with the input \(u_k\)
\begin{align}
    \begin{split}
        x_{k+1} &= f(x_k,u_k) \\
        y_k &= g(x_k)
    \end{split}
\end{align}
will be approximated by a \textit{linear} system
\begin{equation}
    \label{eq:LTI_conti}
    \begin{alignedat}{2}
        z_{k+1} &= A z_k + B u_k &\\
        y_k &= Cz_k &\\
        &&\text{for } z_0 = \phi(x_0),
    \end{alignedat}
\end{equation}
where \(z_k\) is a \textit{lifted} state vector at time-step \(k\).
The nonlinear state vector \(x_k\) will be considered as the output \(y_k\),
so \(y_k := x_k\).
The relationship between the two systems is shown in Fig.~\ref{fig:koopman}.

Having the matrices $A$ and $C$, the matrix $B$ can be optimized
over the whole trajectory, allowing for multiple-step prediction.
The optimization problem can be formulated as
\begin{equation}
    \label{eq:B_min}
    \min \sum_{j=1}^{N_\mathrm{ T }} \sum_{k=1}^{K}
    || g(x_k^j) - \hat{y}_k(x_0^j) ||_2^2,
\end{equation}
where \(N_T\) is the number of trajectories, \(K\) is the number of samples in each trajectory and 
\begin{align}
    \begin{split}
        \label{eq:yhat_k_full}
        \hat{y}_k(x_0^j) = C A^k z_0^j + \sum_{i=0}^{k-1} C A^{k-i-1} B u_i^j,& \\
        \text{for }&z_0^j = \mathbf{\phi}(x_0^j)
    \end{split}
\end{align}
is a prediction of the output vector by the matrices $A$, \(B\) and $C$.

For optimization over shorter window instead of the whole trajectory,
see \cite{cibulka_master_thesis}.

The problem \eqref{eq:B_min} can be also solved as a least-squares problem,
see \cite{Korda_opt_g} for further details.

\subsection{Algorithm summary}
The uncontrolled dynamics is identified first according to Sec.~\ref{sec:uncontrolled} using an uncontrolled dataset.
Then the control is added via the \(B\) matrix, using the approach described in Sec.~\ref{sec:koop_controlled} with a
controlled dataset.
This results in a system
\begin{align}
    \begin{split}
        z_{k+1} &= Az_k + Bu_k \\
        y_k &= Cz_k
    \end{split}
\end{align}
which will be used for linear MPC design.

\section{Identification results}
\label{sec:ident_results}
\subsection{Uncontrolled}
\begin{figure}[htp]
    \centering
    \includegraphics[width=8.8cm]{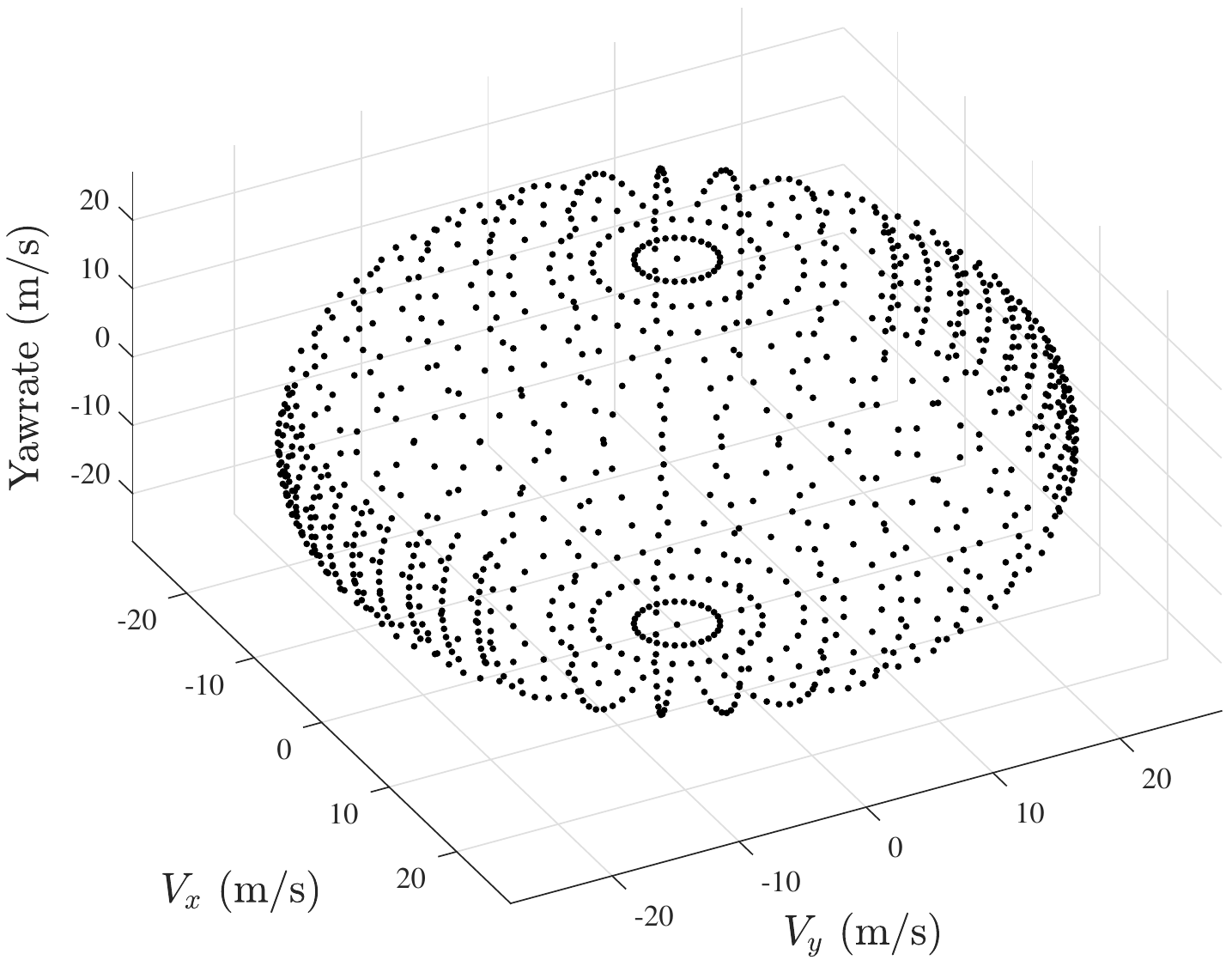}
    \caption{
        The set of initial conditions for the trajectories used for identification.
        The points from this set have a constant kinetic energy \(\ui{E}{k} = \SI{500}{\kilo\joule}\).
    }
    \label{fig:gamma_surface}
\end{figure}

\begin{figure}[htp]
    \centering
    \includegraphics[width=8.8cm]{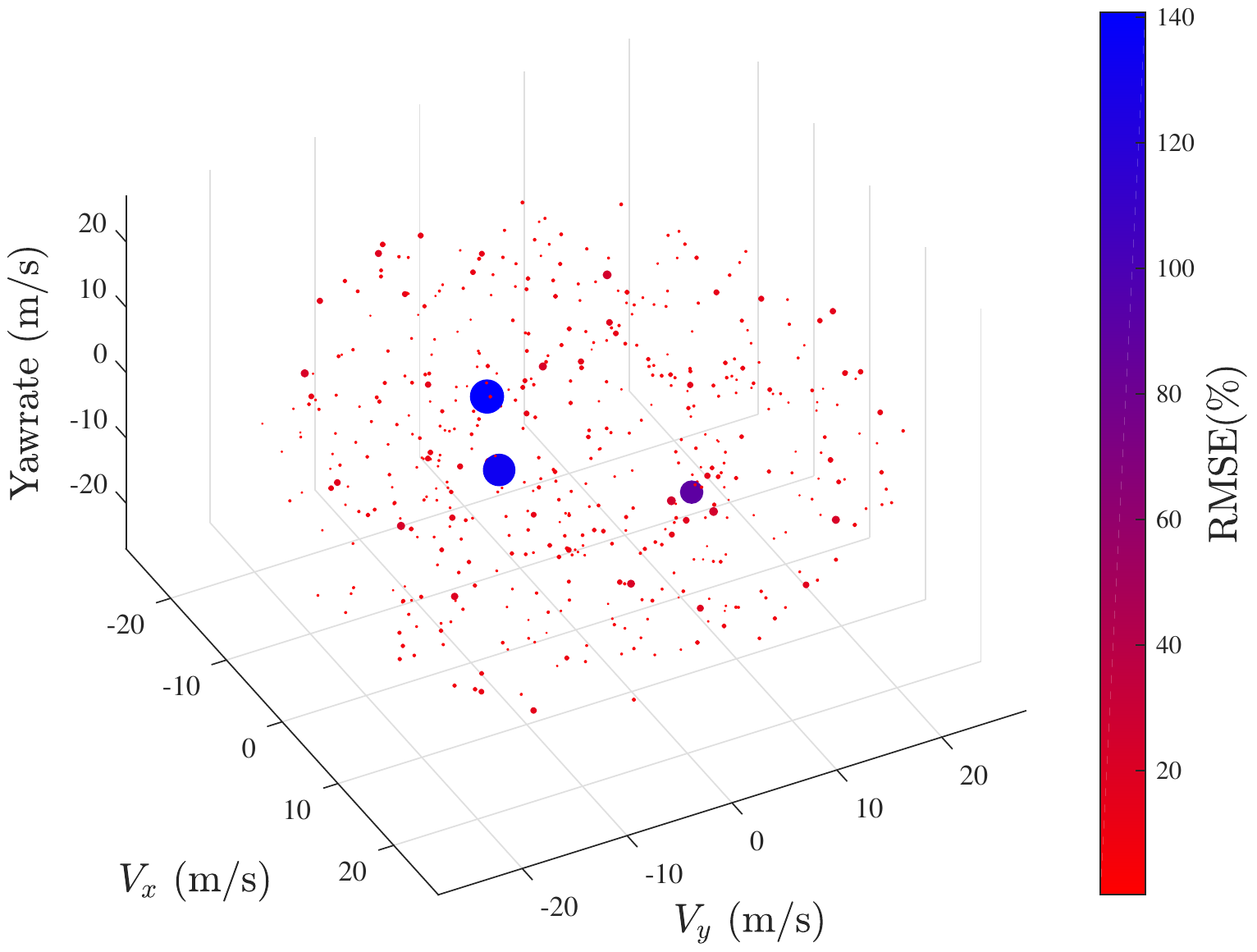}
    \caption{
    Errors of the Koopman operator in the uncontrolled case.
    Each point in the figure corresponds to an initial condition of a \(\SI{0.5}{\second}\) long trajectory.
        The size and color of the points correspond to the prediction error of the associated trajectory. The mean RMSE is \(6\%\) with a
        standard deviation of \(10.2\%\)}
    \label{fig:koop_color_ball}
\end{figure}

\begin{figure}[htp]
    \centering
    \includegraphics[width=8.8cm]{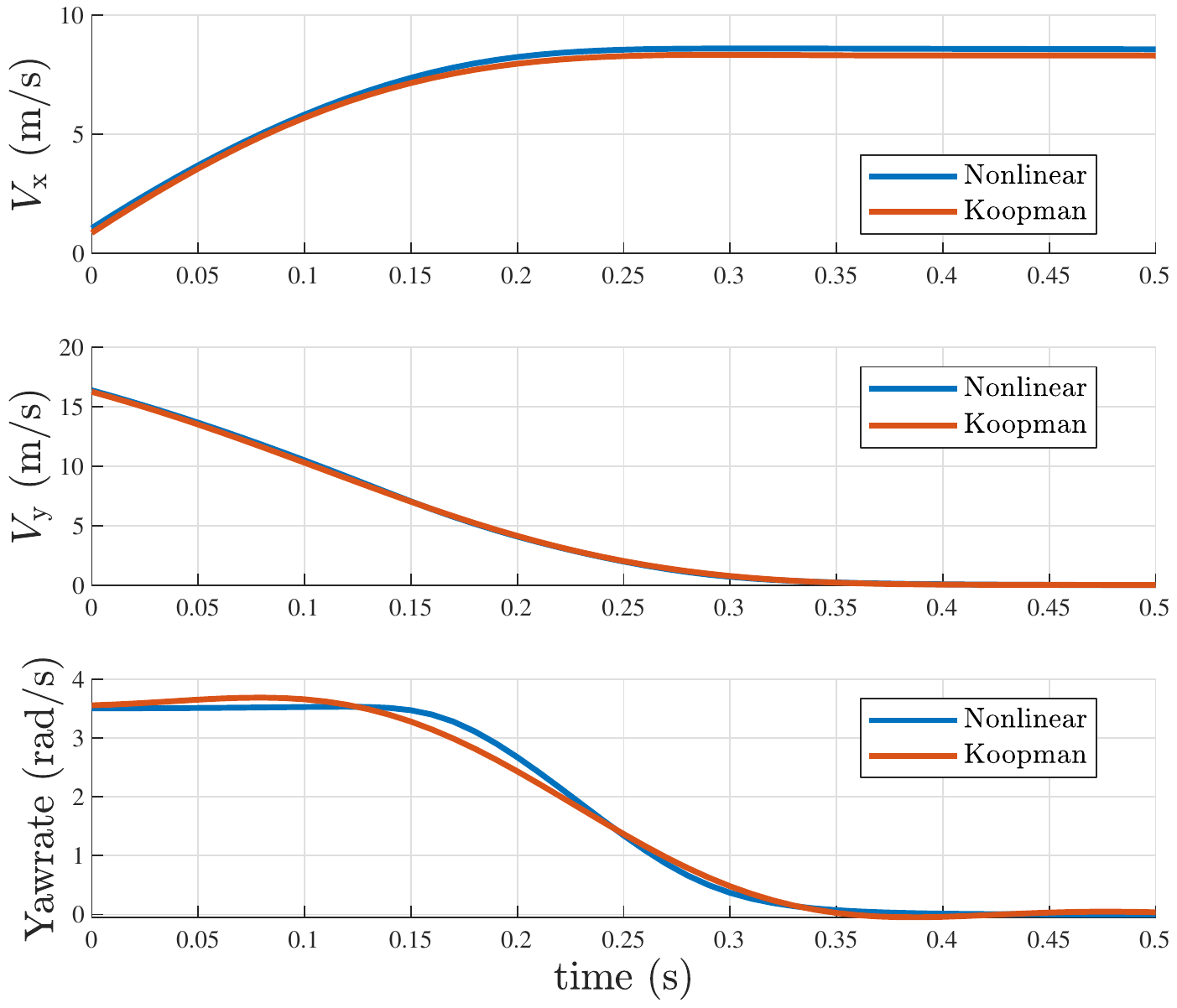}
    \caption{A comparison of the nonlinear and linear system on a trajectory with \(\text{RMSE} = 6\%\) which is equal to the mean RMSE of the whole dataset.}
    \label{fig:koop_trajecotries_mean}
\end{figure}

The model described in Sec.~\ref{sec:singletrack_definition} was
discretized with time-step \(\rm{T_s}=0.01s\)  and approximated by a Koopman operator using the following parameters:
\(N_\Lambda = 51\), \(N_\mathrm{ T } = 1078\) and \(\zeta = 10^{-12}\). The values of the parameters were adopted from \cite{cibulka_master_thesis}, they were 
chosen to provide a sufficient prediction accuracy while keeping computer resource usage at manageable levels.
The starting points for the trajectories were selected from a set with constant kinetic energy \(\ui{E}{k} = \SI{500}{\kilo\joule}\), an
equivalent of a \(\SI{1300}{\kilo\gram}\) car driving straight at \(\SI{100}{\kilo\meter\per\hour}\). The \(\Gamma\) set can be seen in Fig.~\ref{fig:gamma_surface}. Areas
with large \(|v_\mathrm{ y }|\) and low \(|v_\mathrm{ x }|\) (car sliding sideways) contain more points because the vehicle leaves this area of state-space rather quickly,
resulting in sparse data coverage. See \cite{cibulka_koop_identification} for more details.
Results of the uncontrolled dynamics identification can be seen in Fig.~\ref{fig:koop_color_ball} and Fig.~\ref{fig:koop_trajecotries_mean}.
The initial points used for evaluation were randomly generated inside the \(\Gamma\) surface depicted in Fig.~\ref{fig:gamma_surface} and the 
length of the trajectories used for uncontrolled identification was \(0.5~s\), which was the time after which the vehicle model managed to stabilize itself.
Note that this time is rather short because the model defined in Sec.~\ref{sec:singletrack_definition} uses longitudinal slip ratios as inputs and they were set to \(0\)
during the uncontrolled identification. This allowed the tire to generate maximum force in the \(y\) direction which resulted in such short times.
Please see \cite{Pacejka2012} for more information on the tire model.
The starting points with \(||x_0||_2^2 < 8.3\) (\(\SI{8.3}{\meter\per\second}\) \(\dot{=}\) \(\SI{30}{\kilo\meter\per\hour}\)) were rejected from the
testing dataset because the tire model \cite{Pacejka2012} is ill-defined at low speeds.


\subsection{Controlled case}
\begin{figure}[htp]
    \centering
    \includegraphics[width=8.8cm]{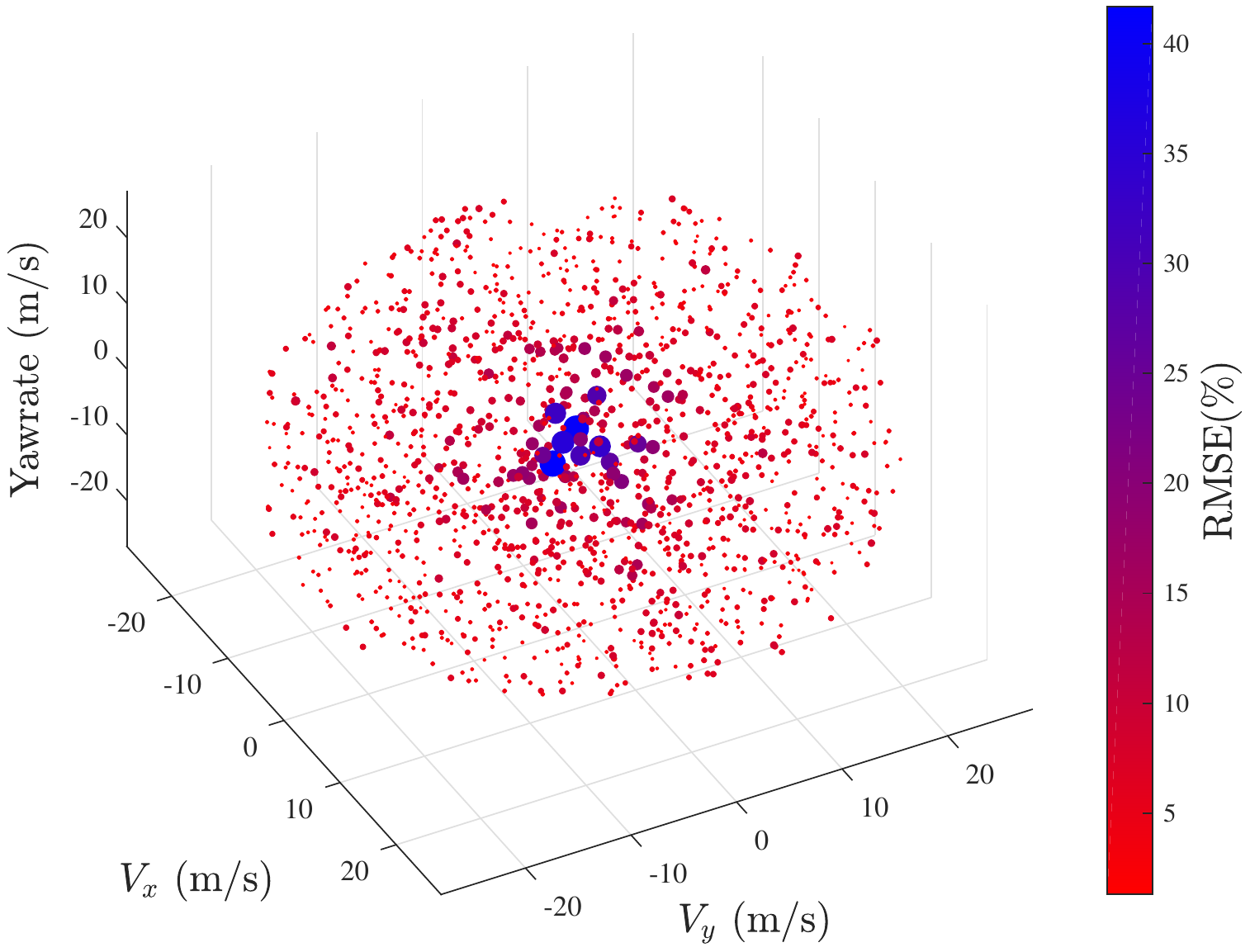}
    \caption{
    Errors of the Koopman operator in the controlled case.
    Each point in the figure corresponds to an initial condition of a \(\SI{0.1}{\second}\) long trajectory with random control inputs.
        The size and color of the points correspond to the prediction error of the associated trajectory. The mean RMSE is \(4\%\) with a
        standard deviation of \(2.7\%\)}
    \label{fig:koop_trajecotries_mean_ctrl}
\end{figure}

Controlled trajectories were generated with randomly generated inputs drawn from a uniform distribution, where
\(\lambda_\mathrm{ r } \in [-1,1]\) and \(\delta_\mathrm{ f } \in [\SI{-30}{\degree},\SI{30}{\degree}]\).
The control horizon for the MPC was chosen as \(\SI{0.1}{\second}\) (adopted from \cite{cibulka_master_thesis}) so the 
matrix \(B\) was optimized on \(\SI{0.1}{\second}\) long trajectories. The mean RMSE was \(4\%\) (the uncontrolled RMSE was \(2.3\%\)).
The distribution of the error can be seen in Fig.~\ref{fig:koop_trajecotries_mean_ctrl}.

\section{MPC}
\begin{figure}[htp]
    \centering
    \includegraphics[width=8.8cm]{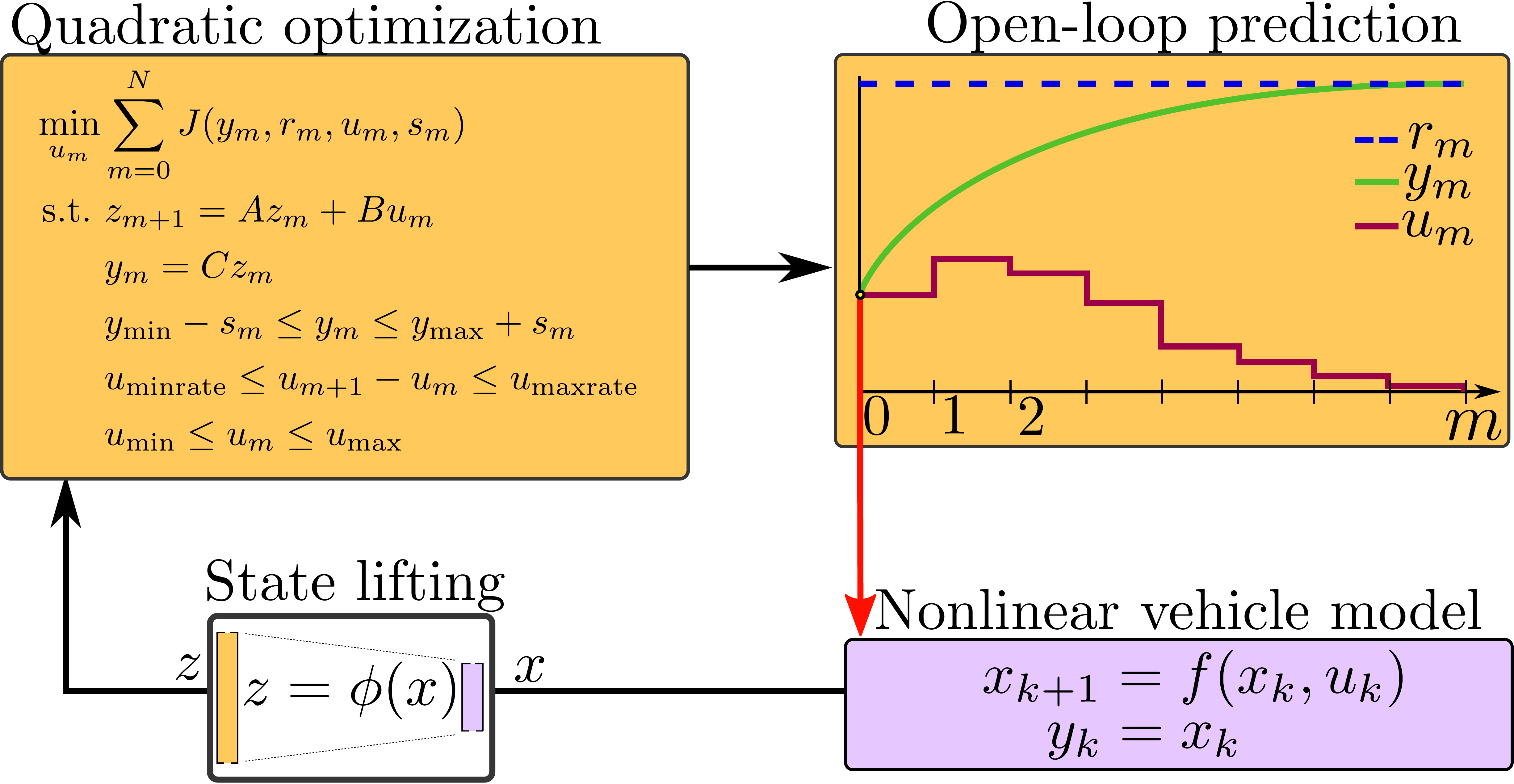}
    \caption{Scheme describing the Koopman MPC algorithm.
        Areas operating in the lifted state-space are depicted in orange color,
        The non-linear space is depicted in violet.
    }
    \label{fig:koop_mpc_scheme}
\end{figure}

The identified system described in Sec.~\ref{sec:ident_results} was used for MPC design.
The MPC based on the Koopman operator will be called Koopman MPC (term first used in \cite{Korda_koop}).
The Koopman MPC framework is depicted in Fig.~\ref{fig:koop_mpc_scheme}.
The Koopman MPC will be compared with MPC based on a locally linearized model, which will be called Linear MPC.
Both MPC regulators are defined as a quadratic optimization problem

\begin{equation}
    \begin{alignedat}{2}
        \label{eq:mpc_definition}
        & \min_{u_m}  \sum_{m=0}^{ N } [(y_m - r_m)^\top \ui{Q}{y} (y_m - r_m) + &u_m^\top R u_m + s_m^\top S s_m] \\
        &\mathrm{s .t.} \\
        &z_{m+1} = Az_m + B u_m                                                     & m=0..N-1&\\
        &y_m = C z_m                                                                & m=0..N-1&\\
        &y_{\mathrm{min}} - s_m \leq y_m \leq y_{\mathrm{max}} + s_m                & m=0..N-1&\\
        & u_{\mathrm{min rate}} \leq u_{m+1} - u_{m} \leq u_{\mathrm{max rate}}     & m=0..N-1&\\
        &u_{\mathrm{min}}  \leq u_m \leq u_{\mathrm{max}}                           & m=0..N-1&,
    \end{alignedat}
\end{equation}

where $Q_\mathrm{ y }$,$S$ and $R$ are positive semidefinite cost matrices, $N$ is the prediction horizon,
$y_{\mathrm{min/max}}$ are soft constraints on the output vector $y_k$ with slack variables $s$
and $u_{\mathrm{minrate/maxrate}}$ are constraints on the system input rates. The only difference between 
Koopman MPC and Linear MPC are the state matrices \(A,B\) and \(C\).
Both MPC regulators were parametrized as follows:
\begin{equation}
    Q_\mathrm{ y } = \begin{bmatrix}
        1 &   &   \\
          & 1 &   \\
          &   & 1
    \end{bmatrix},
    R = \begin{bmatrix}
        0 &     &    &   \\
          & 100 &    &   \\
          &     & 30 &   \\
          &     &    & 0
    \end{bmatrix},
    S = 10^5 \cdot \begin{bmatrix}
        1 &   &   \\
          & 1 &   \\
          &   & 1
    \end{bmatrix}
\end{equation}
\begin{alignat}{2}
    y_{\mathrm{min}}      & = -\begin{bmatrix}
        25 \\2\\2
    \end{bmatrix},
    y_{\mathrm{max}}      &                                & = \begin{bmatrix}
        25 \\2\\2
    \end{bmatrix}, \\
    u_{\mathrm{min}}      & = -\begin{bmatrix}
        0 \\1\\0.45\\0
    \end{bmatrix},
    u_{\mathrm{max}}      &                                & = \begin{bmatrix}
        0 \\1\\0.45\\0
    \end{bmatrix}, \\
    u_{\mathrm{min rate}} & = -\begin{bmatrix}
        0 \\0.1\\0.8\\0
    \end{bmatrix},
    u_{\mathrm{max rate}} &                                & = \begin{bmatrix}
        0 \\0.1\\0.8\\0
    \end{bmatrix}.
\end{alignat}

The scheme of the Koopman MPC is depicted in Fig.~\ref{fig:koop_mpc_scheme}.
The implementation of \eqref{eq:mpc_definition} was done in YALMIP \cite{yalmip}.

\section{Results}
\begin{figure}[htp]
    \centering
        \centering
        \includegraphics[height=12cm]{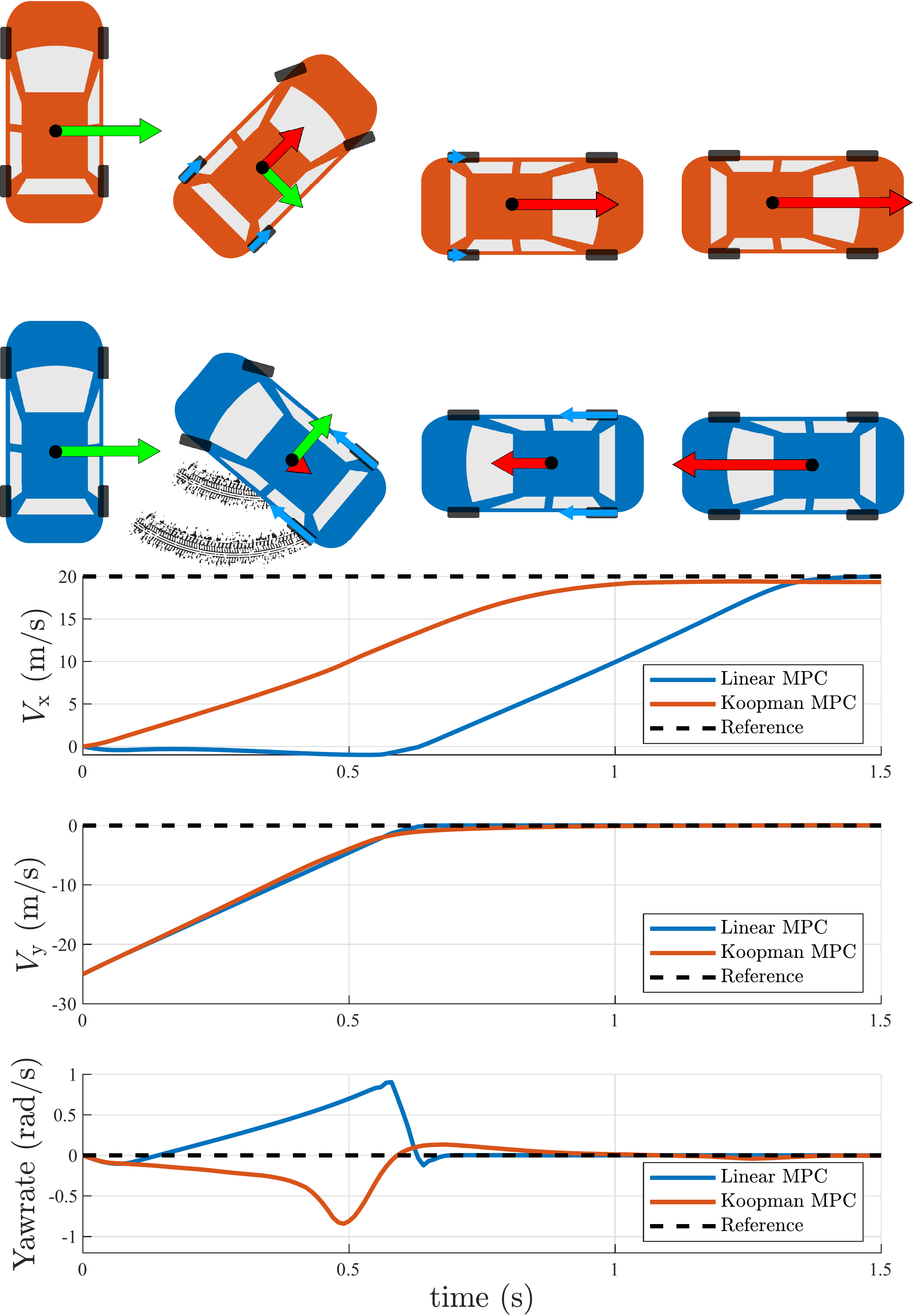}
    \caption{Comparison of maneuvers for recovery from unusual vehicle motion.
        The Koopman MPC stabilized the vehicle faster in one continuous motion.
        The Linear MPC brought the vehicle to full stop and then
        accelerated to reach the desired velocity.
        Notice how both algorithms steered the vehicle in different directions.}
    \label{fig:unusual_attitude}
\end{figure}

It can be seen in Fig.~\ref{fig:unusual_attitude} that the Koopman-controlled vehicle was able to recover from a state where the vehicle
drifts sideways in one continuous motion, unlike the MPC based on local linearization. Notice how each algorithm steered the vehicle in a different direction.
The Koopman MPC steered left in order to shift the momentum from \(\mathrm{ y-axis }\) to \(\mathrm{ x-axis }\) while the locally linearized MPC steered to the right because it was trimmed in state \(
\vect{x}_{\mathrm{0}} = \begin{bmatrix}
    16.7 & 0 & 0
\end{bmatrix}^\top\).
Steering to the right decreases \(v_\mathrm{ y }\) (or increases it in negative direction) in this state.

\begin{figure}[htp]
        \centering
        \includegraphics[height=12cm]{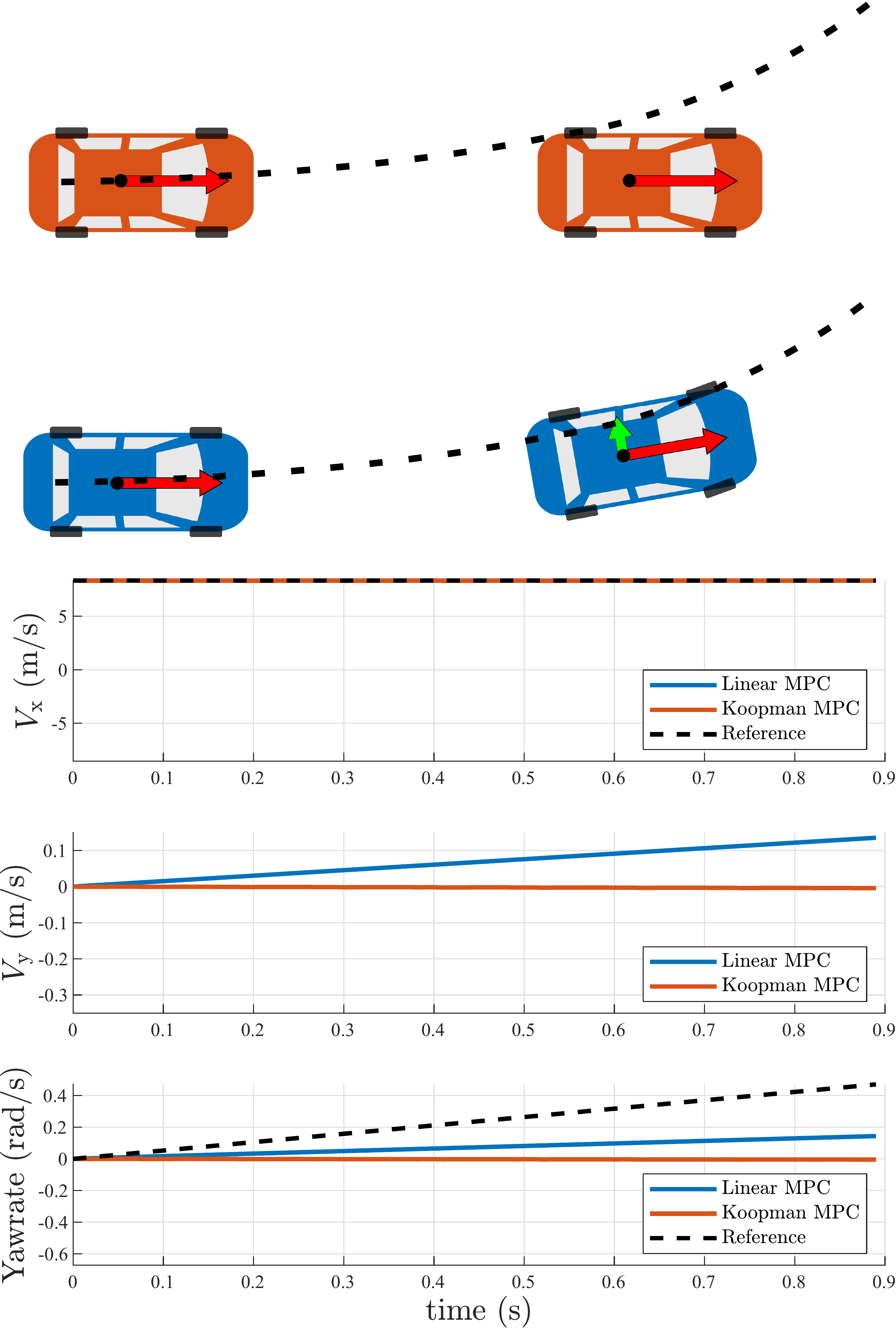}
    \caption{The reference here is increasing yawrate and constant forward velocity \(\ui{v}{x}\). The state \(\ui{v}{y}\) is without any reference.
        The Koopman MPC is unable to track the reference and
        the Linear MPC is performing much better in this case, although it still isn't able to track the reference signal.}
    \label{fig:circle_koop_nefunguje}
\end{figure}

Unfortunately, the Koopman MPC does not always outperform the local linearization. See Fig.~\ref{fig:circle_koop_nefunguje} for example.
In this case, the goal was to steadily increase yawrate while keeping \(v_\mathrm{ x }\) stable. in other words, the vehicle should be driving in an
increasingly tighter spiral while keeping its forward velocity \(v_\mathrm{x}\) the same.

\section{Conclusion}
The Koopman MPC showed very promising results by stabilizing a vehicle from a 90-degree drift while also
preserving energy by shifting the vehicle's already present sideways momentum into a forward momentum.
This result is in stark contrast with the fact the same controller was unable to perform rather simple
steering maneuver.The reason behind this behaviour will be examined in our future work.



\begin{ack}
This research was supported by the Czech
Science Foundation (GACR) under contracts
No. 19-16772S, 19-18424S, 20-11626Y, and 
by the Grant Agency of the Czech Technical 
University in Prague, \\grant No. SGS19/174/OHK3/3T/13.
\end{ack}

\bibliography{ref}             

\begin{thebibliography}{9}
\providecommand{\natexlab}[1]{#1}
\providecommand{\url}[1]{\texttt{#1}}
\providecommand{\urlprefix}{URL }
\expandafter\ifx\csname urlstyle\endcsname\relax
  \providecommand{\doi}[1]{doi:\discretionary{}{}{}#1}\else
  \providecommand{\doi}{doi:\discretionary{}{}{}\begingroup
  \urlstyle{rm}\Url}\fi

\bibitem[{Cibulka et~al.(2019)Cibulka, Hanis, and
  Hromcik}]{cibulka_koop_identification}
Cibulka, V., Hanis, T., and Hromcik, M. (2019).
\newblock Data-driven identification of vehicle dynamics using {K}oopman
  operator.
\newblock \urlprefix\url{http://arxiv.org/abs/1903.06103v1}.

\bibitem[{Cibulka(2019)}]{cibulka_master_thesis}
Cibulka, V. (2019).
\newblock \emph{MPC {B}ased {C}ontrol {A}lgorithms for {V}ehicle {C}ontrol}.
\newblock Master's thesis, CTU in Prague.

\bibitem[{Koopman(1931)}]{Koopman1931}
Koopman, B.O. (1931).
\newblock Hamiltonian {S}ystems and {T}ransformation in {H}ilbert {S}pace.
\newblock \emph{Proceedings of the National Academy of Sciences}, 17(5),
  315--318.
\newblock \doi{10.1073/pnas.17.5.315}.

\bibitem[{Korda and Mezić(2018)}]{Korda_koop}
Korda, M. and Mezić, I. (2018).
\newblock Linear predictors for nonlinear dynamical systems: {K}oopman operator
  meets model predictive control.
\newblock \emph{Automatica}, 93, 149--160.
\newblock \doi{10.1016/j.automatica.2018.03.046}.

\bibitem[{Korda and Mezić(2019)}]{Korda_opt_g}
Korda, M. and Mezić, I. (2019).
\newblock Optimal construction of {K}oopman eigenfunctions for prediction and
  control.
\newblock \urlprefix\url{http://arxiv.org/abs/1810.08733v2}.

\bibitem[{Löfberg(2019)}]{yalmip}
Löfberg, J. (2019).
\newblock {YALMIP}.
\newblock https://yalmip.github.io/.

\bibitem[{Mezi{\'{c}}(2005)}]{Mezic2005}
Mezi{\'{c}}, I. (2005).
\newblock Spectral {P}roperties of {D}ynamical {S}ystems, {M}odel {R}eduction
  and {D}ecompositions.
\newblock \emph{Nonlinear Dynamics}, 41(1-3), 309--325.
\newblock \doi{10.1007/s11071-005-2824-x}.

\bibitem[{Mezi{\'{c}} and Banaszuk(2004)}]{Mezic_control_suggestion}
Mezi{\'{c}}, I. and Banaszuk, A. (2004).
\newblock Comparison of systems with complex behavior.
\newblock \emph{Physica D: Nonlinear Phenomena}, 197(1-2), 101--133.
\newblock \doi{10.1016/j.physd.2004.06.015}.

\bibitem[{Pacejka(2012)}]{Pacejka2012}
Pacejka, H. (2012).
\newblock \emph{Tire and {V}ehicle {D}ynamics}.
\newblock Elsevier LTD, Oxford.
\newblock
  \urlprefix\url{https://www.ebook.de/de/product/18341528/hans_pacejka_tire_and_vehicle_dynamics.html}.

\end{thebibliography}

\appendix
\end{document}